\documentclass[12pt,english]{amsart}
\usepackage[T1]{fontenc}
\usepackage[latin9]{inputenc}
\usepackage{geometry}
\usepackage{color}
\usepackage{babel}
\usepackage{verbatim}
\usepackage{amsthm}
\usepackage{amstext}
\usepackage{amssymb}
\usepackage{esint}
\usepackage[unicode=true]
{hyperref}
\usepackage{breakurl}

\makeatletter
\numberwithin{equation}{section}
\numberwithin{figure}{section}

\@ifundefined{definecolor}
 {\usepackage{color}}{}

\usepackage{babel}

\makeatother

\begin{document}
\title{A brief survey of Nigel Kalton's work on interpolation and related topics}

\author{{Michael Cwikel}}

\address{Cwikel: Department of Mathematics, Technion - Israel Institute of
Technology, Haifa 32000, Israel }

\email{mcwikel@math.technion.ac.il }

\author{{Mario Milman}}

\address{Milman: Department of Mathematics, Florida International University,
Miami, FL 33199, and Department of Mathematics, Florida Atlantic University,
Boca Raton, FL 33431, USA }

\email{mario.milman@gmail.com}

\author{{Richard Rochberg}}

\address{Rochberg: Department of Mathematics, Washington University, St. Louis,
MO, 63130, USA }

\email{rr@math.wustl.edu}

\thanks{The first author's work was supported by the Technion V.P.R.\ Fund
and by the Fund for Promotion of Research at the Technion. The second
author's work was partially supported by a grant from the Simons Foundation
(\#207929 to Mario Milman). The third author's work was supported
by the National Science Foundation under Grant No. 1001488.}
\begin{abstract}
This is the third of a series of papers surveying some small part
of the remarkable work of our friend and colleague Nigel Kalton. We
have written it as part of a tribute to his memory. It does not contain
new results. This time, rather than concentrating on one particular
paper, we attempt to give a general overview of Nigel's many contributions
to the theory of interpolation of Banach spaces, and also, significantly,
quasi-Banach spaces.
\end{abstract}
\maketitle

The heart of interpolation theory is the problem of constructing \textit{interpolation
spaces}. That is, given two Banach (or quasi-Banach) spaces, $A_{0}$
and $A_{1}$ (compatibly contained of course in some larger space),
how can you construct and describe interpolation spaces, $A_{\ast},$
for the pair $\left(A_{0},A_{1}\right)$? Such spaces $A_{\ast}$
should have the interpolation property that a linear operator, $T,$
which is bounded on $A_{i}$ for $i=0,1$ is automatically bounded
on this intermediate space $A_{\ast}$. Although that classical question
has evolved into a broad multifaceted topic, it remains unified by
a few basic themes. Several of those themes were major foci of Nigel
Kalton's research and it is no surprise that he made interesting and
important contributions to interpolation theory.

We will not give an overview of interpolation theory. Basic books
on the topic include \cite{BS}, \cite{BL}, \cite{BK}, \cite{KPS},
\cite{O}, \cite{T}, as well as some earlier books, for instance\ \cite{BB},
and \cite{P}, where one can see viewpoints which guided development
of the general theory. Furthermore, of course, there is the excellent
introduction \cite{iob} to the topic written by Nigel%
\footnote{The fact that in most places in this survey and also in our previous
surveys \cite{CwikelMMilmanMRochbergR2014,CwikelMMilmanMRochbergR2014B}
we allow ourselves the informality of addressing Nigel Kalton by his
first name should be understood as part of our expression of our great
admiration of his work, of our warm friendship for him, and our deep
regret that he is no longer with us.%
} himself, with Stephen Montgomery-Smith, which emphasizes the role
of interpolation theory in the study of the geometry of Banach spaces,
a topic of continuing interest to him.

Nigel's body of work is not a collection of isolated results. It is
a richly interconnected web, and tracing individual strands cannot
capture its richness. That being said, we will now describe several
themes which cover part of Nigel's work on interpolation.

\section{\label{sec:Theory}Theory}

Some of Nigel's contributions were direct additions to the structural
framework of interpolation theory.

Many natural questions at the foundation of interpolation theory are
now well understood, but several remain recalcitrant. Nigel made important
contributions to two of them. One question concerns propagation of
compactness; if the operator $T,$ mentioned earlier, is not only
bounded on $A_{i},$ $i=0,1,$ but also compact on $A_{0};$ must
it then be compact on the interpolation space $A_{\ast}?$ The answer
depends on the details of the method used to construct $A_{\ast}$
For one of the basic classical methods, the complex method, the question
remains unanswered. That question, a favorite of one of the current
authors, is studied in their joint paper \cite{ioc} which we hope
to briefly discuss in a future survey in this series. Another fundamental
question, studied in \cite{cco} asks what conditions on $(A_{0},A_{1})$
insure that the classical real method of interpolation is adequate
to construct all possible interpolation spaces associated with the
pair $(A_{0},A_{1})$. It is classical that the answer is yes for
the pair $(L^{1}.L^{\infty})$ \cite{c2}. (A different way of answering
a closely related question can be found in \cite{MityaginB1965}.)\texttt{
}In \cite{cco} Nigel takes our understanding much further for the
situation in which $L^{1}$ is replaced by a more general rearrangement
invariant space.

The paper \cite{au} is joint work with the authors of this survey.
We have already discussed it in \cite{CwikelMMilmanMRochbergR2014B}.
It gives an approach to interpolation which unifies the real and complex
methods of interpolation (historically the two most important) as
well as several of the less well known variants, in a way which enables
commutator estimates in the spirit of those studied in \cite{do}
to be obtained in a unified way.\textbf{ }(See earlier work by Svante
Janson \cite{JansonS1981} for two other ways of unifying a wide range
of interpolation methods. We also mention that Gilles Godefroy has
discussed the work of \cite{do} (and also several other parts of
Nigel's work) in \cite{Godf07} and \cite{GodfG2012}, and we have
discussed it in \cite{CwikelMMilmanMRochbergR2014}.) 

Questions involving the functorial properties of interpolation constructions
can be quite subtle and results about them can be very useful. For
instance, if a certain construction produces an interpolation space
$A_{\ast}$ from the spaces $A_{i},i=0,1,$ and if the $A_{i}$ have
subspaces $B_{i},i=0,1$, will it be true that the same construction
produces an interpolation space $B_{\ast}$ from the $B_{i}^{\prime}s$
which has the ``expected, natural'' relation to $A_{\ast}$? The
general question was studied by Svante Janson in \cite{J}. In particular
cases, finding a positive answer to this general question can be exactly
the tool needed to solve some other problem. That was done by Nigel
in \cite{cio} and then again, jointly with Sergei Ivanov in \cite{ios}.

As more and more techniques and results of classical analysis were
extended to quasi-Banach spaces, it became important to have a set
of interpolation tools which could work in these new contexts. For
instance, in the mid part of the last century, the theory of singular
integral operators focused on the operators as acting on the Lebesgue
spaces $L^{p}\left(\mathbb{R}^{n}\right),$ $1<p<\infty.$ However,
by the latter part of the century the scope of these studies had expanded
to include viewing these operators as acting on the Hardy spaces,
$H^{p}\left(\mathbb{R}^{n}\right)$, $0<p<\infty;$ a scale of spaces
which coincides with the Lebesgue scale for $p>1$, which is stable
under interpolation, and which, for $p<1,$ consists of quasi-Banach
spaces. One of the tools in the extension of the classical results
to the range $0<p\leq1$ was the theory of complex interpolation of
quasi-Banach spaces, including versions of this theory developed by
Nigel and several of his coauthors. We will have more to day about
this below, in Section \ref{sec:QuasiBanach}. 

Of course Nigel was a leading expert on quasi-Banach spaces (see \cite{qb})
and also studied many other aspects of them which we will not mention
here.

\section{$H^{\infty}$ Functional Calculus, Sectorial Operators}

Before the systematic development of approximation theory or semi-group
theory or interpolation theory; there was the basic question of how
to quantify the degree of smoothness of a function. Counting the number
of derivatives which a function has gives a rough, integer based,
scale; adequate for many purposes. However both curiosity and necessity
rapidly lead to consideration of possible alternatives to, or refinements
of, this scale. There are various intuitively attractive ways to do
this, but difficulties arise when filling in the technical details,
when comparing various approaches to each other, and when exploring
the boundaries of applicability of these techniques.

Consider a toy model question; how can we quantify the degree of smoothness
of a function in $L^{2}=L^{2}\left(\mathbb{T}\text{,}d\theta\right).$
\ We could say that a function in $L^{2}$ is ``smooth'' if the
differentiation operator $D=-i\frac{d}{d\theta}$ maps the function
into $L^{2}.$ The operator $D$ is not bounded but it is closed and
densely defined; it is a plausible starting point for a theory. If
we want to define a fractional order of smoothness using fractional
powers of $D$ we notice that $D$ is not a positive operator and
hence it is not clear how we would want to define, say, $D^{1/2}$.
A standard next step is to instead work with an associated positive
operator of comparable ``size'', for instance $\left\vert D\right\vert =\left(D^{\ast}D\right)^{1/2}$
or $\left(I+D^{\ast}D\right)^{1/2}$ which has the additional property,
at times very convenient, of being bounded below. These operators
are ones for which the classical functional calculus for unbounded
operators on Hilbert space works extremely well and we can use that
calculus to define the fractional powers. We can then use membership
in the domains of definition of the various operators $\left\vert D\right\vert ^{\alpha},$
$\alpha>0$ to define a smoothness scale with a continuous parameter.
We can easily also prove expected results such as: If $0<\alpha<\beta$
and if $g$ has smoothness $\beta$ then $\left\vert D\right\vert ^{\alpha}g$
has smoothness $\beta-\alpha.$

This analysis can be taken quite far; but it is deceptively easy because
passing to Fourier coefficients gives a convenient and explicit diagonalization
of $D$ (and therefore also of $\left|D\right|$) and the smoothness
scale is explicitly defined via simple conditions on Fourier coefficients.
However, even while considering this example, it becomes clear that
more substantial theoretical tools would be needed in more general
situations. What to do if the operator of primary interest is not
positive but just has, in some sense, positive real part? In general,
what functional calculus can be used to define the fractional powers
and show they have the desired computational properties? What can
be done if we measure the size of a function with, say, a Banach space
norm rather than a Hilbert space norm? Do the domains of fractional
powers form an interpolation scale? How do fractional powers of the
operator act on domains of other fractional powers? \ Is it possible
to perturb from a detailed analysis of a simple operator and obtain
results for a modification of that operator?

The history of research in this area includes work to develop holomorphic
functional calculi for various classes of operators on Banach spaces
(to, for instance, define fractional powers), work on the theory of
operator semigroups (again, think of the fractional powers) and the
theory of interpolation of Banach spaces (how are the domains of the
fractional powers related?). Recent overviews of the area include
\cite{Ha} and \cite{KW}. Many of the most substantial developments
in the area were restricted to operators on Hilbert space, for instance
\cite{ANM}. This was in part because, quoting from Nigel's joint
paper \cite{ow} with Tamara Kucherenko, 
\begin{quote}
``\textit{It is an important observation that the theory of sectorial
operators on a Hilbert is, in general, simpler and more easily applicable
than in general Banach spaces. This is mainly due to the characterization
of Hilbert spaces as certain interpolation spaces related to an operator
with an $H^{\infty}$ functional calculus.}'' 
\end{quote}
As the quote suggests, one of Nigel's particular interests in this
area was in expanding the theory beyond the restrictions to Hilbert
spaces. He did that and more. His contributions include \cite{ar},
\cite{ow}, and \cite{pa}, An overview of some of his work in the
area is in \cite{soa}.

\section{\label{sec:QuasiBanach}Quasi-Banach Spaces}

Many interesting and useful classes of functions studied in analysis
are normed spaces and even Banach spaces. However there are also naturally
occurring complete topological vector spaces whose topologies are
defined by a functional $\left\Vert \cdot\right\Vert $ which has
the usual properties of a norm, except that it does not satisfy the
triangle inequality. 
\begin{equation}
\left\Vert x+y\right\Vert \leq\left\Vert x\right\Vert +\left\Vert y\right\Vert .\label{eq:triangle}
\end{equation}
A basic example is, for $p$, $0<p<1,$ the Lebesgue space $L^{p}\left(S,d\mu\right)$
of $\mu$ measurable functions $f$ for which $\int_{S}\left\vert f\right\vert ^{p}d\mu<\infty$.
The functional $\left\Vert f\right\Vert _{p}=\left(\int_{S}\left\vert f\right\vert ^{p}d\mu\right)^{1/p},$
does not satisfy (\ref{eq:triangle}), the space is not a Banach space,
the Hahn-Banach theorem does not hold, and it is not clear if the
space has a non-trivial dual space (the answer depends on $\mu).$
Nevertheless $\left\Vert \cdot\right\Vert _{p}^{p}$ does satisfy
the triangle inequality 
\begin{equation}
\left\Vert x+y\right\Vert _{p}^{p}\leq\left\Vert x\right\Vert _{p}^{p}+\left\Vert y\right\Vert _{p}^{p}\label{eq:pnorm}
\end{equation}
and hence can be used to define a metric, $d(x,y)=\left\Vert x-y\right\Vert _{p}^{p},$
and $L^{p}$ is a complete metric space and is an example of a quasi-Banach
space. The Hardy spaces, $H^{p},$ for $0<p<1,$ are another family
of examples of quasi-Banach spaces.

As we already indicated at the end of Section \ref{sec:Theory}, the
basic questions of interpolation theory continue to be natural and
interesting for quasi-Banach spaces. However one of the fundamental
approaches to interpolation, the complex method of interpolation \cite{C},
does not extend comfortably to quasi-Banach spaces. The reason is
that at the center of the proofs of the fundamental interpolation
inequality is an application of the following inequality, which can
be deduced from the maximum principle for Banach space valued holomorphic
functions: If $f$ is a function which maps the closed unit disk $\mathbb{\bar{D}}$
continuously into a Banach space $X$ and is holomorphic in the interior
of the disk, then 
\begin{equation}
\left\Vert f(0)\right\Vert _{X}\le\int_{0}^{2\pi}\left\Vert f(e^{it})\right\Vert _{X}\frac{dt}{2\pi}.\label{eq:maxz}
\end{equation}
This crucial fact fails for quasi-Banach spaces%
\footnote{It was already known, relatively early in the history of this topic,
that there exist quasi-Banach spaces for which the maximum principle
itself fails, even if we allow a multiplicative constant greater than
$1$ in the relevant inequality. One explicit example of such a space
was given by Aleksandrov \cite[Proposition 4.2, p. 49]{AleksandrovA1981}
(Note that Aleksandrov also connects a corollary of this result with
some earlier work \cite{KaltonN1977} of Nigel.) Another non-explicit
existence proof was obtained, apparently at about the same time, by
the anonymous referee of \cite{PeetreJ1982} and presented in Section
4 (pp,~260--261) of that paper. The proof there is sufficiently ingenious
to make us wonder if that referee could have been Nigel. Who knows?
In the light of these results, it is quite intriguing (though apparently
not helpful for complex interpolation) that, as shown by Nigel in
Section 5 of \cite{pf}, \textit{every} quasi-Banach space satisfies
a slight variant of (\ref{eq:max3}), where instead of integration
over a circle, the integration is performed over an annulus (which
can be arbitrarily thin).%
}. The quasi-norms of holomorphic functions%
\footnote{It is not immediately obvious what is the most appropriate definition
of a holomorphic quasi-Banach space valued function. Here we are implicitly
using the definition via power series adopted by Nigel and authors
of some earlier papers. (See e.g.~the first page of \cite{af} for
this definition and the reason for using it.) But there is at least
one other reasonable alternative definition: When each element of
the quasi-Banach space $X$ is itself a complex valued function on
some underlying measure space $\Omega$, then a given function $f:\mathbb{D}\to X$
of the complex variable $z$ can be considered to be holomorphic (in
what might be called the ``pointwise'' sense) if, at almost every
$\omega\in\Omega$, the value $f(z,\omega)\in\mathbb{C}$ of that
function at $\omega$ depends holomorphically on $z$. Some results
about complex interpolation in the context of this alternative ``pointwise''
definition can be seen, for example, in \cite{GrafakosLMastyloM2014}
and the references therein. %
} which take values in such a space $X$ need not satisfy (\ref{eq:maxz}).
They can also fail to satisfy a less stringent variant of (\ref{eq:maxz}),
namely 
\begin{equation}
\left\Vert f(0)\right\Vert _{X}\le C\int_{0}^{2\pi}\left\Vert f(e^{it})\right\Vert _{X}\frac{dt}{2\pi}\label{eq:max3}
\end{equation}
for some constant $C=C(X)$.

Nigel addressed these matters in an impressive series of papers, some
of the later ones written with co-authors including Loukas Grafakos,
Svitlana Mayboroda and Marius Mitrea. In the initial papers of this
series he built on previous work of a number of mathematicians to
carry out an in-depth study of holomorphic functions taking values
in quasi-Banach spaces, (See \cite{af} and \cite{pf}. Cf.~also
the work of Philippe Turpin \cite{TurpinP1974}.) This enabled him
to then proceed and develop a coherent theory of complex interpolation
of these kinds of spaces. The paper \cite{sro} was his first work
in this direction. As its title indicates, Nigel already dealt here
with interpolation, not only of pairs of quasi-Banach spaces, but
also, more generally, of \textit{families} of such spaces (following
\cite{TabaccoA1988}), and much of his subsequent work would also
apply to families as well as pairs. He identified the condition (\ref{eq:max3})
as the fundamental tool that was needed for developing complex interpolation
for quasi-Banach spaces (and also studied it for other purposes).
He had earlier (in Theorems 3.7 and 4.1 on pp.~305--306 of \cite{pf})
characterized the quasi-Banach spaces for which (\ref{eq:max3}) holds -- their 
function defining the quasi-norm must be equivalent to a
plurisubharmonic function. He and his co-authors also showed that
many of the quasi-Banach spaces of classical analysis such as Hardy,
Besov, Sobolev, and Triebel-Lizorkin spaces satisfy this condition.
(See Section 9 of \cite{KMMioh}.) In particular, detailed accounts
of Nigel and his co-author's version of the theory of complex interpolation
can be found in Section 3 of \cite{KMsri} and in Section 7 of \cite{KMMioh}.
Once the basic tools of the theory were in place, many applications
to classical analysis and partial differential equations became possible.
Nigel and several of his coauthors presented some of these in \cite{srm},
\cite{KMsri}, \cite{KMMioh} and \cite{mc}. For example, \cite{KMMioh}
deals extensively with applications to PDE. 

Nigel's and his coauthors' general and systematic approach to complex
interpolation of quasi-Banach spaces was preceded by quite a number
of earlier papers by many other mathematicians (also including the
authors of this survey) proposing various ways of dealing with this
challenging topic. Some of these were restricted to particular kinds
of spaces, such as $H^{p}$ spaces. The history of this topic is extensive
and complicated and we apologize for not attempting to do it justice
here. We can at least mention that (i) on page 3913 of \cite{KMsri}
Nigel and Marius briefly describe some of the approaches used in some
of these earlier papers, and (ii) on page 157 of \cite{KMMioh} Nigel,
Svitlana and Marius give several references to work where other versions
of complex interpolation have been used for dealing with quasi-normed
Hardy spaces and $BMO$.

\section{Traces and Commutators\label{sec:TracesAndCommutators}}

As the title of this survey indicates, we also wish to briefly discuss
some of Nigel's work in some other topics having some connection with
interpolation theory. We shall do that in this section. Nigel worked
with traces and commutators for most of his career, from his Memoir
\cite{nc} in 1988 to \cite{to} which he did not live to see in print.
His early work in the area was intertwined with his work on interpolation
theory, especially the Memoir and the paper \cite{do}. (As already
mentioned, surveys of the results of \cite{do} can be found in \cite{Godf07,GodfG2012,CwikelMMilmanMRochbergR2014}.)
In his fundamental paper on traces and commutators, \cite{tc}, Nigel
takes the ideas of those previous works much further. However, being
aware that the paper would be of interest to an audience unfamiliar
with interpolation theory, he does not emphasize his path through
the earlier work.

We will not here review Nigel's extensive work on traces and commutators
which includes \cite{tc}, \cite{nc}, \cite{c}, \cite{sc}, \cite{to},
\cite{fs}, and \cite{ri}. However we will point out a particular
set of results which have far flung resonances which we find very
intriguing. So, we will tell that story, but extremely informally.
Some more or less implicit hints about its connections with interpolation
theory may be found in Section 8 of \cite{do} and in the introductory
section of \cite{tc}. 

In the mid 20th Century the University of Chicago was a center of
research in commutative harmonic analysis; it was the home of Antoni
Zygmund and Alberto Calder\'on, and, for a while, Eli Stein and Charles
Fefferman. By the 1950's and 60's much of the classical theory of
Fourier analysis and been extended from the $L^{2}$ spaces of the
line and circle to the corresponding $L^{p}$ spaces, $1<p<\infty$,
and also to functions of several variables. However there were some
fundamental limitations to the theory as it existed then. One problem
was the lack of a satisfying endpoint theory at $p=1$. It was known
that in many ways the classical theory failed to extend to the associated
$L^{1}$ spaces. For functions of one variable, $n=1,$ much of the
theory did go through if, instead of working on $L^{1},$ one worked
with the Hardy space $H^{1},$ the subspace of $L^{1}$ formed by
boundary values of functions holomorphic in the upper halfplane. With
this insight came new questions, what was the appropriate analog of
the Hardy space for functions of $n$ variables and what were the
underlying real variable notions that had made the boundary traces
of holomorphy so surprisingly useful in one dimension? The theory
for the case where the upper halfplane is replaced by the upper halfspace
of $\mathbb{R}^{n+1}$ was thoroughly developed, notably in \cite{FeffermanCSteinE1972}.
Other work in the years that followed led to the introduction of the
``real-variable, atomic Hardy space'', a subspace of functions in
$L^{1}$ which carry the local cancellation properties at the heart
of the classical theory. This viewpoint is presented in \cite{CW};
\cite{S} gives a more recent picture.

Using the general constructions of \cite{CW}, one builds inside of
$L^{1}(\mathbb{R}^{1})$ the subspace ``atomic $H^{1}$'', which
can be denoted by $H_{\mathtt{at}}^{1}=H_{\mathtt{at}}^{1}(\mathbb{R})$
and which is a natural $p=1$ endpoint for many results. Those constructions
are quite general, they work for any metric measure space with certain
resemblances to finite dimensional Euclidean space. In particular,
associated with $\ell^{1}=\ell^{1}\left(\mathbb{Z}_{\mathbf{\geq0}}\right)$
there is the atomic Hardy $h_{\mathtt{at}}^{1}.$ Furthermore, associated
with that space is a symmetrized version, $h_{\mathtt{sym}}^{1}$
which is, roughly, the space of sequences with rearrangements in $h_{\mathtt{at}}^{1}.$
Precisely, $h_{\mathtt{sym}}^{1}$ is the space of sequences $a=\left(a_{1}.a_{2},...\right)$
such that, with $\left(\tilde{a}_{1}.\tilde{a}_{2},...\right)$ denoting
the same sequence of numbers but rearranged so that $\left\vert \tilde{a}_{n}\right\vert $
is nonincreasing, we have
\[
\left\Vert a\right\Vert =\sum\left\vert a_{n}\right\vert +\sum\left\vert \frac{\tilde{a}_{1}+\tilde{a}_{2}+...+\tilde{a}_{n}}{n}\right\vert <\infty,
\]
 In fact $h_{\mathtt{sym}}^{1}$ is a quasi-Banach space with quasi-norm
equivalent to $\left\Vert \cdot\right\Vert .$ This space, which is
studied in \cite{tc}, \cite{do}, and \cite{nc} is, as is discussed
in those references, a space which can be naturally viewed as a discrete
analog of the space of rearrangements of functions in $H_{\mathtt{at}}^{1}(\mathbb{R}).$
Hence, in contexts where one is measuring the size of sequences using
rearrangement invariant functionals, $h_{\mathtt{sym}}^{1}$ is a
candidate to be a useful extension to $p=1$ of the scale of spaces
$\ell^{p}\left(\mathbb{Z}_{\mathbf{\geq0}}\right),$ $p>1.$

The Schatten classes $\mathcal{S}_{p}$ are the Banach spaces of linear
operators on a Hilbert space with the property that their singular
values are in the sequence space $\ell^{p}$, $0<p<\infty$. It was
a question of interest in operator theory to know how the Schatten
classes, and other analogously defined operator ideals were interrelated
under commutation. For instance it was known that, for $p>1,$ we
have $\mathcal{S}_{p}=\left[\mathcal{S}_{2p},\mathcal{S}_{2p}\right]$
where, in the notation of that area $\left[A,B\right]$ is the closed
linear span of the set commutators $\left[a,b\right]$, $a\in A,b\in B,$
and this was known that this fails at $p=1.$ For background on this
see \cite{dfww} and \cite{W}.

Recalling that the Schatten classes are, in some sense, analogs of
the Lebesgue classes, it is perhaps not surprising that a Hardy type
space, in this case, $h_{\mathtt{sym}}^{1}$ ,leads to the natural
analog for $p=1$. Let $\mathcal{S}_{h_{\mathtt{sym}}^{1}}$ be the
space of compact operators on a Hilbert space with the property that
their sequence of eigenvalues $\left\{ \lambda_{n}\right\} \in h_{\mathtt{sym}}^{1}.$
(Note that $h_{\mathtt{sym}}^{1}$ is rearrangement invariant and
hence $\mathcal{S}_{h_{\mathtt{sym}}^{1}}$ is unitarily invariant.
However, in contrast to $\ell^{p}$ spaces, $h_{\mathtt{sym}}^{1}$
is not a solid space. With this in mind the definition of $\mathcal{S}_{h_{\mathtt{sym}}^{1}}$
is in terms of eigenvalues, not singular values ( = approximation
numbers.)) In fact in \cite{tc} Nigel proves that $\mathcal{S}_{h_{\mathtt{sym}}^{1}}=\left[\mathcal{S}_{2},\mathcal{S}_{2}\right].$
That is 
\[
A\in\mathcal{S}_{h_{\mathtt{sym}}^{1}}\text{ if and only if }A=\sum\left[B_{i},C_{i}\right];\text{ }B_{i},C_{i}\in\mathcal{S}_{2},
\]

One reason for pointing out this result is to note the analogy with
a different result in the theory of Hardy spaces. Now let $H^{p},$
$0<p<\infty$ be the classical Hardy spaces associated to the unit
ball $\mathbb{B}^{n}$ in $\mathbb{C}^{n}.$ That is, functions in
$L^{p}$ of the surface of $\mathbb{B}^{n}$ which are boundary values
of functions holomorphic in $\mathbb{B}^{n}.$ The following theorem
is from \cite{CRW}
\[
A\in H^{1}\text{ if and only if }A=\sum B_{i}C_{i};\text{ }B_{i},C_{i}\in H^{2}.
\]
The visual analogy between these two results is clear. Some hints
at the deeper relations are seen by comparing Nigel's papers which
were mentioned with the ideas of \cite{CLMS}. It seems possible that
there are more systematic interrelationships to be found.

\section{Finally}

As we already indicated at the outset, it is of course impossible,
in the limited framework of this survey, to come anywhere near doing
full justice to all of Nigel's broad and extraordinary contribution
to interpolation theory. Ideally, many more of Nigel's papers, including
ones which we could only mention briefly here, would have been described
in more detail in some parallel documents. Among our regrets we note,
for example, that the papers \cite{oa} and \cite{ci}, no less fine
than those we have mentioned, did not fit comfortably into our narrative.
And, given the extent of Nigel's total research output and the interplay
of different topics within it, it seems almost certain that we have
overlooked other papers of his which are relevant in one way or another
to this survey.

\end{document}